\documentclass[a4paper,twoside,11pt]{amsart}
\usepackage{newtxmath}
\usepackage{amssymb,amsmath}

\usepackage{type1cm}        

\usepackage{mathtools}

\def\int{\mathrm{int}}

\makeindex       

  \title{A propos des travaux d'Albert Lautman}

   \author{Bernard Teissier}
      
\begin{document}
\maketitle
Lorsque j'ai d\'ecouvert les \'ecrits d'Albert Lautman dans \cite{L}, j'ai imm\'ediatement \'et\'e tr\`es int\'eress\'e car cela me semblait bien plus proche de ce que je ressentais comme (jeune) math\'ematicien que tout ce que j'avais lu jusqu'alors dans le domaine de la philosophie des Math\'ematiques. Cette id\'ee qu'il fallait s'int\'eresser \`a ce qui structurait les math\'ematiques, \`a ce qui faisait leur unit\'e sous-jacente, plut\^ot qu'\`a des questions comme "le mode d'existence des objets math\'ematiques" ou \`a leur seule structure logique. Je crois que j'ai toujours pens\'e \`a la logique, disons \`a la Frege, comme \`a l'\'echafaudage qui emp\^eche les b\^atisseurs d'une cath\'edrale de faire une chute mortelle lorsqu'ils font un pas risqu\'e comme celui que fit Cantor. Il ne faut pas confondre l'\'echafaudage et la cath\'edrale. L'id\'ee de r\'eduire les Math\'ematiques \`a leur structure logique m'a toujours r\'evuls\'e.\par Bien s\^ur la logique math\'ematique moderne est d'une autre nature et en particulier explore des objets math\'ematiques d'un point de vue issu de pr\'eoccupations logiques. Par exemple, comment mieux comprendre la complexit\'e des preuves ou le rapport entre le langage et la syntaxe utilis\'es pour d\'efinir axiomatiquement des objets math\'ematiques et la nature proprement math\'ematique de ces objets.\par
Pour revenir \`a Lautman, j'\'etais heureux de trouver aussi dans ses \'ecrits l'id\'ee d'une dynamique des math\'ematiques dont le moteur est en partie le dialectique des concepts et dont le r\^ole est en partie de rendre manifeste l'unit\'e profonde des Math\'ematiques. D'autres, bien plus comp\'etents que moi, pr\'esenteront la richesse des id\'ees de Lautman, dont je ne per\c cois qu'une petite partie. Ce que je veux faire ici est exposer  ce que je pense \^etre la cause d'une frustration que j'ai ressentie d\`es mes premi\`eres lectures mais dont la cause m'\'echappait alors.\par\noindent
Cete cause est que pour moi, la Math\'ematique est une science humaine. Je veux dire par l\`a que son existence et tout son d\'eveloppement d\'ependent avant tout de notre humanit\'e.  C'est d'ailleurs vrai aussi pour la logique, qui sait r\'ef\'echir sur sa nature m\^eme, en prenant en compte ses rapports avec l'informatique th\'eorique (voir \cite{G}, et une bonne partie de l'oeuvre de J.-Y. Girard, qui introduit sciemment beaucoup d'humanit\'e dans la logique). L'objectivit\'e des Math\'ematiques est r\'eelle en ce sens que les \'enonc\'es vrais resteraient vrais pour toute entit\'e acceptant les pr\'emisses et les r\`egles logiques mais ces \'enonc\'es n'auraient peut-\^etre aucune signification pour cette entit\'e, qui ne les \textit{comprendrait} pas. Une philosophie des Math\'ematiques qui postule un caract\`ere absolu ou "objectif" de leur naissance et de leur d\'eveloppement est \`a mon avis dans une orni\`ere.\par Or, malgr\'e les grands m\'erites du travail de Lautman, je n'y ai pas trouv\'e cette id\'ee. Il se peut qu'elle soit cach\'ee dans des vocables dont je comprends mal le sens philosophique, comme "logique transcendentale" ou "esth\'etique transcendentale" mais lorsqu'il \'ecrit que \textit{l'objectivit\'e des Math\'ematiques r\'eside dans leur participation \`a une r\'ealit\'e plus haute et plus cach\'ee, …un monde des id\'ees}, je ne peux pas le suivre s'il attribue vraiment aux id\'ees une existence objective. \par
Depuis au moins trois d\'ecennies, de nombreux chercheurs ont commenc\'e \`a explorer un nouveau  r\^ole de l'inconscient dans la philosophie des Math\'ematiques. Cela est d\^u en grande partie au d\'eveloppement des sciences cognitives, qui permet de commencer \`a objectiver notre rapport inconscient au monde et de cesser de le dissimuler derri\`ere un rideau de p\'eriphrases comme "intuition sensible" ou "r\'ealit\'e". Ces d\'eveloppements auraient certainememt beaucoup int\'eress\'e Lautman. Un exemple majeur est le texte \cite{P} de Jean Petitot, qui cite au d\'ebut le texte suivant de David Hilbert dans \textit{\"Uber das Unendlich}:\par\noindent
\textit{Kant avait d\'ej\`a pour doctrine que les Math\'ematiques ont un contenu ind\'ependant de la logique et qu'elles ne peuvent \^etre fond\'ees sur la logique seulement; c'est ce qui condamnait d'avance les tentatives de Frege et de Dedekind. La condition pr\'ealable \`a l'usage des inf\'erences logiques est l'existence d'un donn\'e dans la perception….}\par\noindent
Je me place r\'esolument dans cette lign\'ee, non en tant que philosophe des sciences, mais en tant que math\'ematicien int\'eress\'e par des questions de nature plut\^ot philosophique \`a propos des Math\'ematiques.\par
J'ai propos\'e dans \cite{T1} de distinguer, en ce qui concerne les constructions et les r\'esultats math\'ematiques, les fondements de la v\'erit\'e et les fondements de la signification. Le probl\`eme des fondements de la v\'erit\'e a donn\'e naissance \`a de magnifiques d\'eveloppements motiv\'es en grande partie par le souci de rigueur "absolue" donnant la certitude d'\'eviter les contradictions. Mais comme ces d\'eveloppements sont le fruit d'une r\'eflexion humaine, ils ont bien s\^ur leur propre signification, qui d'ailleurs depuis quelques d\'ecennies se rapproche de la g\'eom\'etrie avec la th\'eorie des topos.\par
En ce qui concerne la signification, j'ai propos\'e d'explorer une voie qui n'existait pas du temps de Lautman, et que l'on pourrait appeler la "subjectivit\'e objective". Il s'agit de l'id\'ee que les fondements de la signification se trouvent dans l'\'enorme quantit\'e d'exp\'erience du monde par nous et nos anc\^etres qui est contenue dans notre inconscient et dans la structure de notre cerveau. Et cette exp\'erience du monde est suffisamment stable quand on passe d'un individu \`a l'autre pour que nous puissions \'echanger des "assemblages de signification"  tr\`es complexes. C'est dans cette stabilit\'e, qui est un fait d'exp\'erience, que r\'eside \`a mon sens l'objectivit\'e de ce que nous appelons des id\'ees. \par\noindent
Ce qui m'a conduit dans cette direction est la constatation que comprendre un \'enonc\'e ou sa d\'emonstration est de la nature d'une illumination et non le r\'esultat d'un cheminement logique. Et ce qui m'a encourag\'e \`a continuer est le d\'eveloppement des sciences cognitives. 
Il me semble que nous comprenons une d\'emonstration quand nous avons extrait du texte un tissu d'interpr\'etations en cascade des objets math\'ematiques impliqu\'es en termes de notre exp\'erience primitive du monde qui est compatible avec cette exp\'erience.\par
C'est pendant notre formation de math\'ematicien, et plus tard dans la pratique,  que nous apprenons \`a donner sens \`a des objets math\'ematiques complexes, en termes de notre exp\'erience inconsciente du monde et aussi bien s\^ur de notre exp\'erience des Math\'ematiques elles-m\^emes. C'est cette mani\`ere de donner sens qui fait que nous \textit{comprenons} au sens \'etymologique, des objets tr\`es complexes. C'est \`a cause du lien tr\`es \'etroit entre notre perception du monde, la mani\`ere dont cette perception a model\'e notre cerveau, et la cr\'eation des Math\'ematiques, que je pense qu'il est absurde, en principe, de s'\'etonner comme Wigner de "The unreasonable effectiveness of Mathematics". C'est s'\'etonner que l'\'ecorce colle \`a l'arbre. Mais \'evidemment les processus qui conduisent \`a cette "unreasonable effectiveness" restent extr\^emement myst\'erieux et fascinants. \par
Mon exemple favori, d\'evelopp\'e dans \cite{T1}, est que notre perception inconsciente du monde contient deux d\'efinitions de la droite: la droite vestibulaire, qui correspond \`a un \'etat extremal d'une assembl\'ee de neurones quand nous marchons \`a vitesse constante dans une direction fixe, et la droite visuelle, qui correspond \`a un \'etat extremal d'une autre assembl\'ee de neurones du cortex visuel lorsque notre oeil d\'etecte un segment de droite. Et notre syst\`eme perceptif inconscient \'etablit une correspondance entre les deux, que j'ai appel\'ee "Isomorphisme de Poincar\'e-Berthoz" (voir \cite{B}, \cite{B2}) et qui permet par exemple de repr\'esenter le temps (qui est une mesure provenant de la droire vestibulaire) comme param\'etrant la droite r\'eelle, qui vient de la droite visuelle. Cette fusion perceptive de la droite vestibulaire avec la droite visuelle est l'origine, \`a mon avis, de beaucoup de d\'eveloppements math\'ematiques. Einstein disait qu'une de ses intuitions fondamentales \'etait de se d\'eplacer sur un rayon de lumi\`ere. Je dirais qu'Einstein avait une relation particuli\`erement riche avec son inconscient perceptif. Notre syst\`eme perceptif relie ainsi, \`a sa mani\`ere, le discret et le continu. En vertu du principe de relativit\'e, la progression sur la droite vestibulaire ne peut \^etre mesur\'ee que de deux mani\`eres: une mesure continue par la distance parcourue ou le temps pass\'e, en supposant connue la vitesse (c'est \`a 30 minutes de marche d'ici) et d'autre part le nombre de pas, qui est une mesure discr\`ete. Il me semble voir l\`a l'origine de la question de savoir si le temps, ou le continu, est form\'e d'indivisibles et en m\^eme temps l'origine du concept de trajectoire. En rapport avec cette question des indivisibles, les math\'ematiciens voulaient comprendre ce qui rend le continu diff\'erent d'un ensemble de points "en vrac"; est-ce la structure d'ensemble ordonn\'e (les pas du marcheur ou le temps) et le fait qu'il est un ensemble de fronti\`eres (les coupures de Dedekind) ? Cela est vraiment proche de la proto-pens\'ee et a finalement conduit \`a la th\'eorie des ensembles! \par\noindent Identifier le continuum visuel avec le continuum temporel a des cons\'equences \'enormes, comme le concept de trajectoire param\'etr\'ee par le temps. Et il y a beaucoup plus, rien que dans notre syst\`eme visuel (voir \cite{P-T}, \cite{P2}). Par exemple notre syst\`eme visuel est construit de mani\`ere \`a d\'etecter, dans notre perception de l'espace,  le transport parall\`ele d'Elie Cartan. \par
Par ailleurs, il me semble que la dialectique, concept philosophique bien compris, est loin d'\^etre le seul moteur du d\'eveloppement des Math\'ematiques. Il y a me semble-t-il tout un \'eventail de pulsions primitives, pour la plupart inconscientes, dont j'ignore si elles sont "c\^abl\'ees" dans notre cerveau, mais qui sont \'egalement assez stables quand on passe d'un individu \`a l'autre. Appelons cela la "pens\'ee primitive" ou "proto-pens\'ee"\footnote{Motiv\'e par la "vision de bas niveau" des sp\'ecialistes de la physiologie de la vision, dans \cite{T1} j'avais utilis\'e le vocable "pens\'ee de bas niveau" mais cela avait suscit\'e des erreurs d'interpr\'etation}. Voici des exemples:\par\noindent
\begin{itemize}
\item Comparer des objets comparables: au vu de deux objets, je sais imm\'ediatement lequel est plus gros, lequel est plus loin, \textit{sans me poser la question}.
\item D\'etecter des r\'egularit\'es, d\'etecter des traits analogues dans des objets. D\'etecter la r\'ep\'etition des r\'esultats d'exp\'eriences semblables. Faire mentalement des it\'erations.
\item Capacit\'e de projeter des repr\'esentations mentales, de simplifier et surtout d'abstraire m\^eme en l'absence de langage.
\item Une recherche obstin\'ee des causes, des origines, et plus g\'en\'eralement, une \textit{curiosit\'e} insatiable. Bon nombre de probl\`emes sont du type: si A implique B, est-il vrai que B implique A?
\end{itemize}
Cette liste est loin d'\^etre exhaustive, et l'\'etude de cette proto-pens\'ee serait fascinante. Nous en partageons une partie avec les primates.\par
A nouveau, il me semble que refuser d'admettre le r\^ole de cette proto-pens\'ee dans le d\'eveloppement des Math\'ematiques laisse la r\'eflexion dans une orni\`ere. Par exemple, notre exp\'erience primitive du monde comprend les ombres, et leur forme abstraite est la projection. Il me semble que la "mont\'ee vers l'absolu", que Lautman a eu le m\'erite d'identifier parmi les mouvements importants des Math\'ematiques, fait sens pour nous \`a cause de cela et du d\'esir de simplification. \par
De m\^eme, l'exp\'erience primitive de la marche et la proto-pens\'ee de l'it\'eration font que nous pouvons sans trop de difficult\'e donner un sens au concept d'infini m\^eme si ses propri\'et\'es peuvent faire d\'ebat et si la question de "l'existence de l'infini" n'a pu \^etre r\'esolue qu'assez tardivement - et d'ailleurs magnifiquement - dans un cadre math\'ematique. \par Si pour nous l'id\'ee de "marcher ind\'efiniment dans la m\^eme direction" fait sens, celle "d'avoir march\'e ind\'efiniment jusqu'ici" n'est pas acceptable aussi facilement (probl\`eme de l'origine) et il me semble qu'il faut chercher l\`a, et dans l'apparition relativement r\'ecente de la soustraction comme op\'eration susceptible d'it\'eration, la source de la signification de la notion d'ensemble bien ordonn\'e.\par
Pour conclure, il me semble que, \'etant donn\'e le d\'eveloppement des sciences cognitives,  la philosophie des Math\'ematiques, si elle se veut aussi proche de la nature m\^eme des Math\'ematiques que celle de Lautman, ne peut manquer de s'int\'eresser au r\^ole que jouent les proto-objets math\'ematiques comme le proto-continuum cr\'e\'e par l'isomorphisme de Poincar\'e-Berthoz ainsi que les proto-pens\'ees dont j'ai \'enum\'er\'e une petite partie. Je pr\'ecise, s'il est n\'ecessaire, qu'il ne s'agit nullement d'une approche r\'eductionniste. Je suis convaincu que les processus en jeu du point de vue physiologique sont d'une complexit\'e qui litt\'eralement d\'epasse l'entendement et que nous ne pouvons avoir qu'une vision assez floue de la mani\`ere dont ils influencent les Math\'ematiques. Mais cette vision serait d\'ej\`a passionnante ! \par\medskip
\textit{Je remercie chaleureusement les organisateurs du colloque de m'avoir permis d'exprimer ces r\'eflexions d'un non-philosophe.}

 Bernard Teissier \par\medskip\noindent
 Universit\'e Paris Cit\'e and Sorbonne Universit\'e, CNRS, IMJ-PRG, F-75013 Paris, France \par\noindent
bernard.teissier@imj-prg.fr
 

\begin{thebibliography}{A}
\bibitem{B} A. Berthoz, Le cerveau, le mouvement, et les espaces, in \textit{Neurosciences cognitives}, sous la direction de Mehdi Khamassi, \'editions De Boeck Sup\'erieur, 2021, 25-55.
\bibitem{B2} A. Berthoz, La marche, le cerveau, et l'espace, les g\'eom\'etries du corps en marche, in \textit{Le g\'enie de la Marche} sous la direction de Sabine Chardonnet Darmaillacq, coll. Colloque de Cerisy, Hermann, 2016, 295-315.
\bibitem{G} J.-Y. Girard, \textit{Le point aveugle (2 volumes)}, Hermann, 2006-2007.
\bibitem{L} A. Lautman, \textit {Essai sur l'unit\'e des Math\'ematiques et divers \'ecrits}, Union G\'en\'erale d'Edition, Paris 1977.
\bibitem{P} J. Petitot, \textit{Refaire le "Tim\'ee". Introduction \`a la philosophie math\'ematique d'Albert Lautman}, Revue d'Histoire des Sciences, XL, 1, 79-115 (1987).
\bibitem{P} J. Petitot, Continu et objectivit\'e, in \textit{Le labyrinthe du continu}, Colloque de Cerisy, J.-M. Salanskis et H. Sinaceur, \'editeurs, Springer, 1992.
\bibitem{P2} J. Petitot, \textit{Neurog\'eom\'etrie de la vision} Les Editions de l'Ecole Polytechnique, 2008.
\bibitem{P-T} J. Petitot et Y. Tondut, \textit{Vers une Neuro-g\'eom\'etrie. Fibrations corticales, structures de contact et contours subjectifs modaux}, Num\'ero sp\'ecial de Math\'ematiques, Informatique et Sciences Humaines, 145, 5-101, EHESS, Paris.
\bibitem{T1} B. Teissier, G\'eom\'etrie et Cognition: l'exemple du continu, in \textit{Ouvrir la logique au monde", actes de l'Ecole th\'ematique CNRS-LIGC " Logique et interaction; vers une g\'eom\'etrie du cognitif}, Cerisy Septembre 2006, dirig\'ee par J.-B. Joinet et S. Tron\c con, Hermann, "Visions des sciences", Paris 2009. Disponible par: https://webusers.imj-prg.fr/~bernard.teissier/documents/Cerisy06final5.pdf
\bibitem{T2} B. Teissier, Mathematics and Narrative: why are stories and proofs interesting?, in \textit{Circles disturbed: the interplay of Mathematics and Narrative}, edited and introduced by Apostolos Doxiadis and Barry Mazur, Princeton University Press, 2012. Disponible par: https://webusers.imj-prg.fr/~bernard.teissier/documents/Teissier.corr.doc
 \end{thebibliography}
  \end{document}